\documentclass[a4paper,10pt]{article}
\def \Z {{\mathbf {Z}}}
\def \R {{\mathbf {R}}}
\def \N {{\mathbf {N}}}

\def \e {{\bf 1}}
\def\u{\bigsqcup}
\def\eps{\varepsilon}
\title{ Minimal Self-Joinings,  Bounded Constructions, \\ and  Weak Closure of Ergodic   Actions. }
\author{V.V. Ryzhikov}
\date{12.12.2012}

\begin{document}
\large
\maketitle
%\begin{abstract}{\Large 
For a weakly mixing bounded  rank-one  construction it is proved  the disjointness of its powers.  For non-rigid constructions we get  minimal self-joinings.  Examples of non-mixing rank one actions  with explicit weak closure are proposed. 
\section{Introduction} 
Consider rank-one measure-preserving transformations  of a Probability space $(X,\mu)$.
In connection with Thouvenot's question on minimal self-joinings
(MSJ) of mildly mixing rank-one transformations, for a class of transformations (mildly mixing bounded constructions) we prove the property MSJ. This class includes well-known classical and modified  Chacon's transformations (see\cite{JRS}).

The "springs" of our proof are similar to ones of \cite{JRS},  although  we do not consider generic points, and never  use Birkhoff's  ergodic theorem. Weak limits are our tools.  In operator terms
the method is following.  Let $P$ be an indecomposable Markov
operator commuting with a rank-one transformation $T$ ($P$ is  extreme point of  Markov centralizer of $T$).
If $P\neq T^i$ for all $i\in \Z$, then there is a sequence ${n_j}$
such that $T^{n_j}\to aP+(1-a)P'$, where  $a\geq \frac 1 2$,
$P'$ is some Markov operator.  In fact there is a sequence of
sets $Y_j$ such that $\mu(Y_j)\geq \frac 1 2$  and
$$\e_{Y_j}T^{n_j}\to aP, \ \ \  a\geq \frac 1 2.$$
Using non-rigidity of $T$   find a sequence $C_j\subset
Y_j$ for which
$$\e_{C_j}T^{n_j}\to cT^m P,\ \  \ \ c>0,  \ \ m\neq 0.$$
But from  $C_j\subset
Y_j$ we have  $\e_{C_j}T^{n_j}\to cP$ as well. This implies
$$P=T^mP.$$
If $T$ is mildly mixing, then  $T^m$ ($m\neq 0$) is ergodic,
hence, $P=\Theta$ -- the ortho-projector to constant functions in
$L_2(X,\mu)$.
We get the triviality of Markov centralizer. It is the set of all convex sums of the operators
$T^i$ and $\Theta$, consequently, $T$ has  MSJ.

In \cite{JRS} the authors used a similar trick: they 
proved for  a non-trivial self-joining $\nu$ of Chacon's map
the marginal invariance:
$$  \nu =(Id\times T)\nu, $$ 
so $  \nu =\mu\times \mu. $

We give   new  examples of non-mixing rank one actions  with explicit weak closure  that gives minimal self-joinings. 
An action  of such kind appeared in arXiv:1108.0568: for  double staircase transformation $T$  the  weak closure of its powers $Lim(T)$  is $$\{\Theta, 2^{-m}T^n+(1-2^{-m})\Theta \ : m\in\N,\ n\in \Z\}.$$ 
A non-simple description of power weak limits for  Chacon's transformation is given in \cite{JPRR}.

In this note we show that $Lim(T)$ for so-called stochastic Chacon's transformation (see \cite{R}) is (a.s.) 
$$\{\Theta, P^{m}{P^\ast}^{n}T^k \ : m,n\in\N,\ k\in \Z, m+n\neq 0,  P=aI+(1-a)T\}.$$ 

For bounded staircase flow $T_t$  we get
$$Lim(T_t)=
\{\Theta, T_a\prod_{i\in S} P_{m_i}\ :a\in \R,  S\subset \N, |S|<\infty,
 P_m=\int_0^mT_t dt\}.$$ 
\section{Rank-one constructions}
{ A rank-one construction} is determined by an integer $h_1$, a cut-sequence $r_j$ and  a spacer  sequence $\bar s_j$
$$ \bar s_j=(s_j(1), s_j(2),\dots, s_j(r_j-1),s_j(r_j)).$$  We recall its definition.
Let a transformation  $T$ is defined on the step $j$ as  a shift on a collection of disjoint sets (intervals)
$$E_j, TE_j T^2E_j,\dots, T^{h_j}E_j.$$
Cut $E_j$ into $r_j$ sets (subintervals)  of the same measure
$$E_j=E_j^1\u E_j^2\u  E_j^3\u\dots\u E_j^{r_j},$$  
then for all $i=1,2,\dots, r_j$ we  consider so-called columns
$$E_j^i, TE_j^i ,T^2 E_j^i,\dots, T^{h_j}E_j^i.$$
Adding $s_j(i)$ spacers over a $i$-th column we obtain  "column + spacers"
$$E_j^i, TE_j^i T^2 E_j^i,\dots, T^{h_j}E_j^i, T^{h_j+1}E_j^i, T^{h_j+2}E_j^i, \dots, T^{h_j+s_j(i)}E_j^i$$
(all intervals   are disjoint).
For all  $i<r_j$ we "stack"
$$TT^{h_j+s_j(i)}E_j^i = E_j^{i+1}.$$ 
Now one obtains a tower 
$$E_{j+1}, TE_{j+1} T^2 E_{j+1},\dots, T^{h_{j+1}}E_{j+1},$$
where 
 $$E_{j+1}= E^1_j,$$
$$T^{h_{j+1}}E_{j+1}=T^{h_j+s_j(r_j)}E_j^{r_j}, $$
$$h_{j+1}+1=(h_j+1)r_j +\sum_{i=1}^{r_j}s_j(i).$$

\bf Bounded   constructions. \rm  Following \cite {B}  consider  constructions  assuming  that   all $s_j(i)$, $r_j$ are bounded: $0\leq s_j(i)<s$,  $2<r_j<r$ .

(A reader can  slightly generalize the above notion, only requiring   for   all $j$ bounded "derivatives":
 $\max_{1\leq i< r_j} |s_j(i+1)-s_j(i)|<s, \ \ r_j<r.\ )$

If $s_j(r_j)=0$ for   all $j$, we get a "strongly bounded construction", which is  as a rule similar to  modified Chacon's transformation .
\\
\bf Examples. \rm

Odometers. For instance,  $r_j=5,$  $\bar s_j = (2,2,2,2, 0)$ for all $j$.

\rm Classic
Chacon's map:  $r_j=2$, 
  $\bar s_j = (0,1)$ for all $j$. It is not mixing, that is clear because of  the following  weak convergence
$$ 2T^{-h_j}\to \sum_{i=0}^{\infty}\left(\frac T 2 \right)^i. $$

 Modified Chacon's map (\cite{JRS}): $ r_j=3$, 
  $\bar s_j = (0,1, 0)$ for all $j$. Now 
$$ T^{-h_j}\to \frac {I+T} 2.$$

\section{ Joinings} \rm A {\it self-joining} (of order 2) is defined to be a $T\times T$-invariant
measure $\nu$ on $X\times X$ with the marginals  equal to  $\mu$:
$$\nu(A\times X)=\nu(X\times A)=\mu(A).$$
A joining $\nu$ is called ergodic if the dynamical system
$(T\times T, X\times X, \nu)$ is ergodic.
The measures $\Delta^i= (Id\times T^i)\Delta$ 
defined by the formula
$$
	\Delta^i(A\times B) = \mu(A\cap T^iB) 
$$ 
are referred to as {\it off-diagonals measures\/} (for $i\neq 0$). 
If $T$ is ergodic, then $\Delta^i$
are   ergodic self-joinings. 

We say that  $T$  has {\bf minimal self-joinings}
 of   order 2 (and we write $T\in MSJ(2)$)  if $T$  has no  ergodic
joinings except $\mu^{\otimes 2}=\mu\times\mu$ and $\Delta^i$.

The notion of MSJ (of all orders) has been introduced by D. Rudolph \cite{Ru} (see also \cite{T}).  In \cite{JRS} the authors  proved MSJ for 
modified Chacon's automorphism.  It is well-known fact that
MSJ(2)=MSJ for non-mixing maps (Glasner, Host, Rudolph, Ryzhikov).

The property of minimal self-joinings implies {\bf  mild mixing}. We recall that an automorphism $T$ is {\it mildly mixing} if  for any set $A$, $0<\mu(A)<1$,  
$$\limsup_j 
\mu(A\cap T^iA)  < \mu(A) .$$
An automorphism  $T$ is mildly mixing iff it has no rigid factors ( $S$ is rigid, if there is a sequence $m_j\to\infty$
such that $S^{k_j}\to Id$).  

 An automorphism $T$ is  {\it partially mixing}, if
 for for some $\alpha \in (0,1]$ and for all measurable sets $A,B$ 
$$\liminf_i 
\mu(A\cap T^iB)  \geq \alpha \,\mu(A)\mu(B).$$
% \vspace{2mm}
In \cite{KT} the authors proved minimal self-joinings for partially mixing
rank one transformations.
The property of partial mixing implies mildly
mixing. 

J.-P. Thouvenot conjectured  that all 
\bf bounded  mildly mixing  constructions have minimal self-joinings. \rm 
We give a positive answer (generalizing \cite{JRS} and some results from \cite{K}).
\vspace{3mm}

{\bf THEOREM.} \it  Let $T$ be a bounded non-rigid construction.  If it is totally ergodic (all non-zero powers are ergodic), then $T$ has minimal self-joinings.\rm
\vspace{3mm}

\section{Ergodic limits of off-diagonals are trivial}
\it Given  an  ergodic self-joining $\nu$ of a rank one transformation $T$, is   there  a sequence $n_j$ such that for all measurable $A,B$
$$\nu(A\times B) = \lim_{j\to\infty}   \mu(T^{n_j}A\cap B)      \ \ {\bf ?}$$
\rm This  question is due to J. King.

If such a sequence $n_j$ exists,  we call  $\nu$ 
a limit of off-diagonals and write
 $$\Delta_{T^{n_j}}\to \nu.$$

\bf LEMMA      \it Let $T$ be non-rigid totally ergodic bounded construction,  $m_j\to\infty$ and  $\Delta_{T^{m_j}}\to \nu.$  If $\nu$ is ergodic, then $\nu=\mu\times\mu$. \rm

Proof.  We find  $p=p(i)$  such that  $h_p\leq m_j <h_{p +1}$. Now we consider our construction on the step $p$.
Let's remark that $p=p(i)\to\infty$ and $rh_p>h_{p +1}$ for all (large)$p$.  The spacers  over the last column form 
the roof over $(p+1)$-th tower.  The transformation $T$ is non-rigid, hence the roof have to be (asymptotically) non-flat.  This implies the following:   $$T^{m_j}C_j\to cP,$$ and  for some Markov operator $Q$ commuting with $T$ we have 
$$P=(a_0I +a_1T +\dots+ a_{r-1}T^{r-1})Q,$$ where  at list two coefficients, say $a_n$, $a_{n+k}$, $k>0$, have to be both non-zero.   Let us rewrite this in joining terms:
$$\Delta_{T^{m_j}}\to \nu = \sum_{k=0}^{r-1}a_k(Id\times T^k)\eta.$$
 Thus, $\nu$ and $(Id\times T^k)\nu$ are not disjoint.
Assuming $\nu$ to be ergodic   we get
$$  \nu =(Id\times T^k)\nu. $$ 
From the ergodicity of $T^k$ it follows that
$$   \nu =\mu\times\mu $$ 
(in operator terms we say:  $T^kP=P \ \ implies\ \  P=\Theta.$)
Indeed,
$$   \nu(A\times B) =\int_{X\times X} \chi_A\otimes \left(\frac 1 N \sum_{d=1}^N T^{dk}\chi_B\right) d\nu =$$
$$ = \lim_N  \int_{X\times X}\chi_A\otimes \left(\frac 1 N \sum_{d=1}^N T^{dk}\chi_B\right) d\nu=
 \mu(B)\int_{X\times X}\chi_A\otimes 1 d\nu = \mu(A)\mu(B). $$

\section {Joinings as Local Limits}
 \rm Any ergodic self-joining  $\nu$ of a rank one transformation $T$  is  a partial limit of off-diagonals:   there is a sequence $n_j$ such that
$$\Delta_{T^{n_j}}\to \frac {1}{2} \nu  +\frac {1}{2} \nu'.$$
In fact (see  \cite{JRR}), for some 
$\delta \leq \frac {1}{2}$ there is a sequence of sets $Y_j$ in the form 
$$Y_j=\bigcup_{\delta h_j < k < h_j}T^kE_j$$  and a sequence $\{n_j\}$, $n_j \approx (1-\delta)h_j$,   such that 
$$\nu(A\times B) = \lim_{j\to\infty}  \frac {1}{\mu(Y_j)} \mu(T^{n_j}A\cap B\cap Y_j),$$
equivalently, 
$$\frac {1}{\mu(Y_j)} \e_{Y_j}T^{n_j}\ \to\  P.$$
Such limits we call local.
The following lemma shows that sometimes certain local limits become global.
 
\bf LEMMA 1. \it  Let $\e_{Y_j}T^{n_j}\to (1-\delta)P$,  $h_j\leq n_j<h_{j+1}$.
We represent  
$$n_j = q h_{j} +s_j(1)+s_j(2) + \dots + s_j(q) + m_j, \ \  0\leq m_j <h_j. $$

  If $\frac {m_j}{ h_{j}}\to 0$ (or $ \frac {h_{j}-m_j}{ h_{j}}\to 0$),  then 
$T^{m_j}\to P$ (or $T^{m_j-h_j}\to P$). \rm

Proof.  Let $C_j^1$ denote  the first column of $j$-tower.
We see that  $$\mu(C_j^1)>\frac 1 r, \ \ \
T^{n_j}C_j^1 \subset Y_j, \ \ \ \mu(T^{n_j}C^1_j\Delta C^{q+1}_j)\to 0.$$
We get
$$\nu(A\times B) = \lim_{j\to\infty}  \frac {1}{\mu(Y_j)} \mu(T^{n_j}A\cap B\cap Y_j)=$$

$$= \lim_{j\to\infty}  \frac {1}{\mu(C^1_j)} \mu(T^{n_j}A\cap
T^{n_j}C^1_j\cap B)=$$
$$= \lim_{j\to\infty}  \frac {1}{\mu(C^1_j)} \mu(T^{m_j}A\cap
C^{q+1}_j\cap B)
=\lim_{j\to\infty}  \mu(T^{m_j}A
\cap B).$$

\section{On local "breaking"  of non-trivial  self-joinings }

    Assume  $T$ satisfies the conditions of Theorem 1.  Let's find 
sequences of sets  $C_j', C_j''\subset Y_j$ such that 

\ \ \ \ \ \ \  \ \ \ \  $\mu(C_j')\approx \mu(C_j'')\to c>0$,\ \ $  \e_{ C_j'}{}^{{}_\circ} T^{n_j}\to cP$,\ \
$  \e_{ C_j''}{}^{{}_\circ} T^{n_j}\to cP$,
\\ and in addition  for fixed  $ m>0$ and some  sequence $m_j$,  $|m_j|<h_j$,   let
$$  \e_{ C_j'}{}^{{}_\circ} T^{n_j}\approx_w \e_{ C_j'}{}^{{}_\circ}  T^{m_j},\eqno (*) $$
$$   \e_{ C_j''}{}^{{}_\circ} T^{n_j}\approx_w  \e_{ C_j''}{}^{{}_\circ}T^{m_j +m},\eqno (**)$$ 
$$  \e_{  C_j'}{}^{{}_\circ} T^{m_j}\approx_w \e_{ C_j''}{}^{{}_{\circ }} T^{m_j}. \eqno (***)$$
If such $C_j', C_j''$ are found,  then
$$cP\approx_w \e_{ C_j''}{}^{{}_\circ}T^{m_j +m} \approx_w
\e_{ C_j'}{}^{{}_\circ}T^{m_j +m}\approx_w cT^mP, \ \
P=T^mP, \ \ P=\Theta.$$

Repeat the above in joining terms : 
 for all measurable $A,B$
$$\nu(A\times B) =
 \lim_{j\to\infty}   \mu(T^{n_j}A\cap B\cap C''_j)/{\mu(C''_j)}\\
=_{(**)}  \lim_{j\to\infty}  \mu(T^{m_j+m}A\cap B\cap C''_j)/{\mu(C''_j)}=_{(***)}$$
$$  \lim_{j\to\infty}   \mu(T^{m_j}T^m A\cap B\cap C'_j)/{\mu(C'_j)} =_{(*)}
\lim_{j\to\infty}   \mu(T^{n_j}T^m A\cap B\cap C'_j)/{\mu(C'_j)}=\nu(T^m A\times B).$$
Thus, 
$  \nu =(T^m\times Id)\nu,   \ \ \  \nu=\mu\times\mu. $

\bf How to find  $C_j', C_j''$?
\rm
Consider  minimal  $i$ for which $s_j(i)\neq s_j(i+1)$.  Without loss of generality  suppose $i<0.4 r_j$.

{\normalsize In oder to have this, instead of the sequence
$\{j\}$ (of steps)  we may consider  the subsequence  $\{Nj\}$,    where $N$ is fixed.  Now   $$\tilde r_{j}:= \prod _{k=(j-1)N+1}^{jN}r_{k}.$$  
 Let's look at Chacon's map  spacer sequence  $\bar s_j = (0,1, 0)$.
Setting $j:=2j$, we get new spacer sequence  
$\bar{\tilde s}_j= (0,1,0,0,1,1,0,1,1). $
Generally,  given construction $T$,  if for all $N$  we cannot find $i$ s.th.   $s_j(i)\neq s_j(i+1)$, then we see  flat  roofs  of towers,   hence,  $T$ have to be  rigid. It is not our case. Thus, we get 
for some $N$  for all (new) $j$  a desired integer $i=i(j)$  satisfied  $s_j(i)\neq s_j(i+1)$. If $N$ is sufficiently large, we find $i<0.4r_j$. }

If  the image under $T^{n_j}$ of the left half  of $E_j$ ( and the left part of spacers)  is
located far away from the top and bottom of the tower,
then we easy find   $C_j', C_j''$ as below.  Let's remark that in the case of 
the  flat image  at the top of the tower, we are looking for $C_j', C_j''$ at the bottom  (and vice versa, respectively).

 
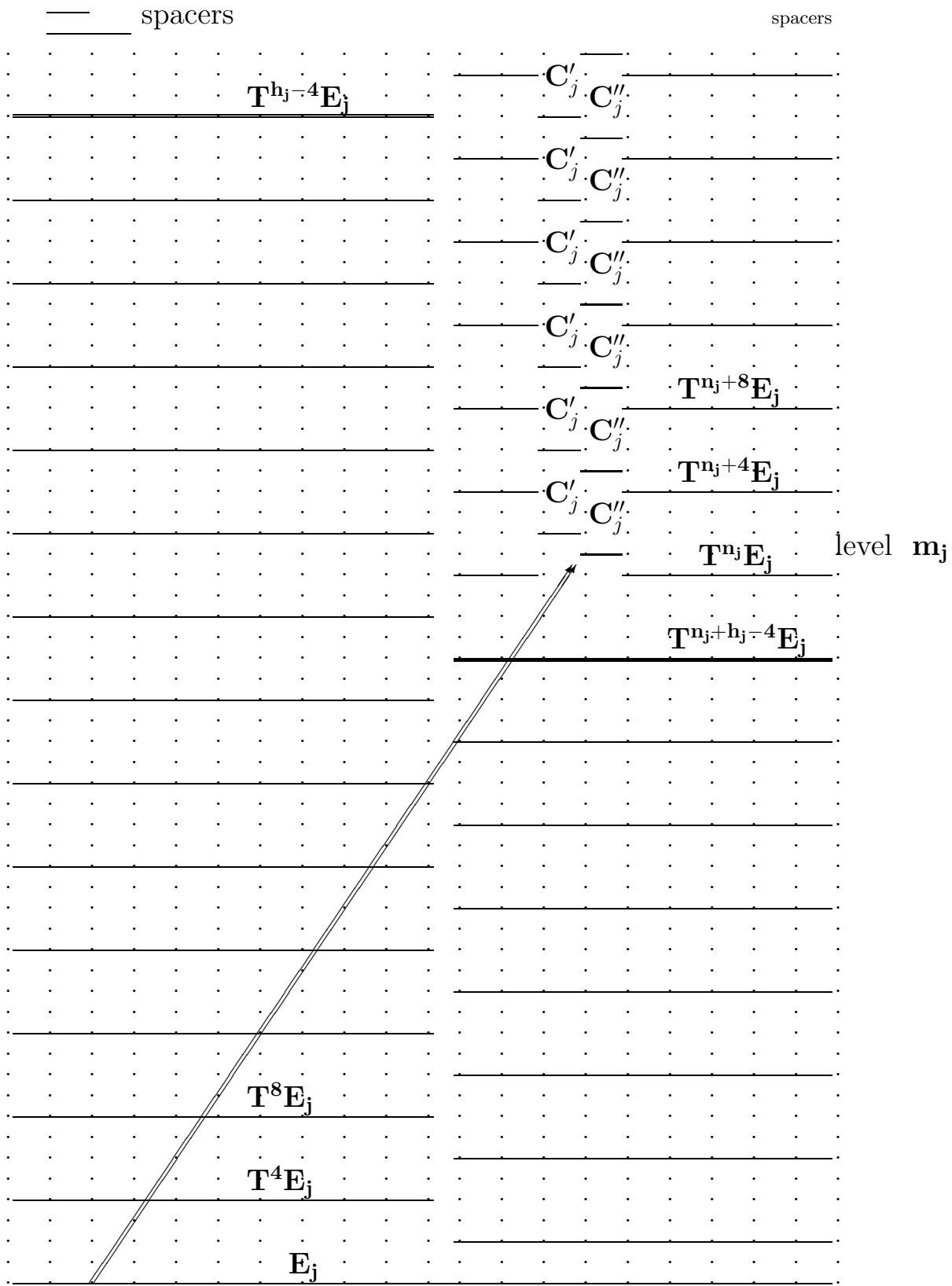
{%\small
\begin{picture}(0, 600)
%\put(200,600){\it \Large The sets $C'_j$, $C''_j$ up }
\put(42,0){\vector(2,3){230}}
\put(0,0){  \line(1, 0){390}}
%\put(-4,600){  \line(1, 0){390}}
\put(16,600){  \line(1, 0){40}}
\put(16,610){  \line(1, 0){20}}
{\Large
\put(65,605){spacers}}

{\small
\put(365,605){spacers}}
%\put(0,0){  \line(0, 1){600}}
%\put(400,0){  \line(0,1){600}}

\multiput(0,0)(20,0){11}{
\multiput(0,0)(0,10){60}%
{${.}$}}
\multiput(215,0)(20,0){10}{
\multiput(0,0)(0,10){60}%
{${.}$}}

%\multiput(400,0)(0,10){60}%
%{${\bf .}$}
\multiput(0,0)(0,40){15}
{  \line(1,0){200}}

\multiput(210,20)(0,40){8}
{  \line(1,0){180}}
\put(210,299){  \line(1, 0){180}}

%\put(250,301){  \line(0, 1){290}}

\put(40,0){\vector(2,3){230}}
\multiput(210,580)(0,-40){7}
{  \line(1,0){40}}
\multiput(270,590)(0,-40){7}
{  \line(1,0){20}}
\multiput(270,590)(0,-40){6}
{  \put(8,-25){{\Large \bf C$''_j$}}}
\multiput(250,560)(0,-40){6}
{  \put(7,15){{\Large \bf C$'_j$}}}
\multiput(250,560)(0,-40){6}
{  \line(1,0){20}}
\multiput(290,580)(0,-40){7}
{  \line(1,0){100}}

\put(0,561){  \line(1, 0){200}}
%\put(0,558){  \line(1, 0){200}}

%\put(210,421){  \line(1, 0){40}}
%\put(210,219){  \line(1, 0){40}}
%\put(290,221){  \line(1, 0){100}}
%\put(290,219){  \line(1, 0){100}}

%\put(250,201){  \line(1, 0){20}}
%\put(250,199){  \line(1, 0){20}}

%\put(270,211){  \line(1, 0){20}}
%\put(270,210){  \line(1, 0){20}}
%\put(270,209){  \line(1, 0){20}}
%\put(130,100){  \li,ne(1,-2){50}}
{\Large
\put(135,5){${\bf  E_j}$}
		\put(115,45){${\bf  T^4E_j}$}
\put(115,85){${\bf T^8 E_j}$}
}

{\Large
\put(330,345){${\bf T^{n_j} E_j}$}
\put(395,350){level $ \ {\bf {m_j}}$}
		\put(320,385){${\bf T^{n_j+4} E_j}$}
\put(320,425){${\bf  T^{n_j+8}E_j}$}

}

{\Large
\put(115,565){${\bf  T^{h_j -4}E_j}$}
		\put(315,305){${\bf  T^{n_j+h_j -4}E_j}$}}
%\put(150,-12){$r_j=19$} }
\end{picture}
%  \begin{center}Fig. 1\end{center}
%%%%%%%%%%%%%%%%%%%%%%%%%%%%%%

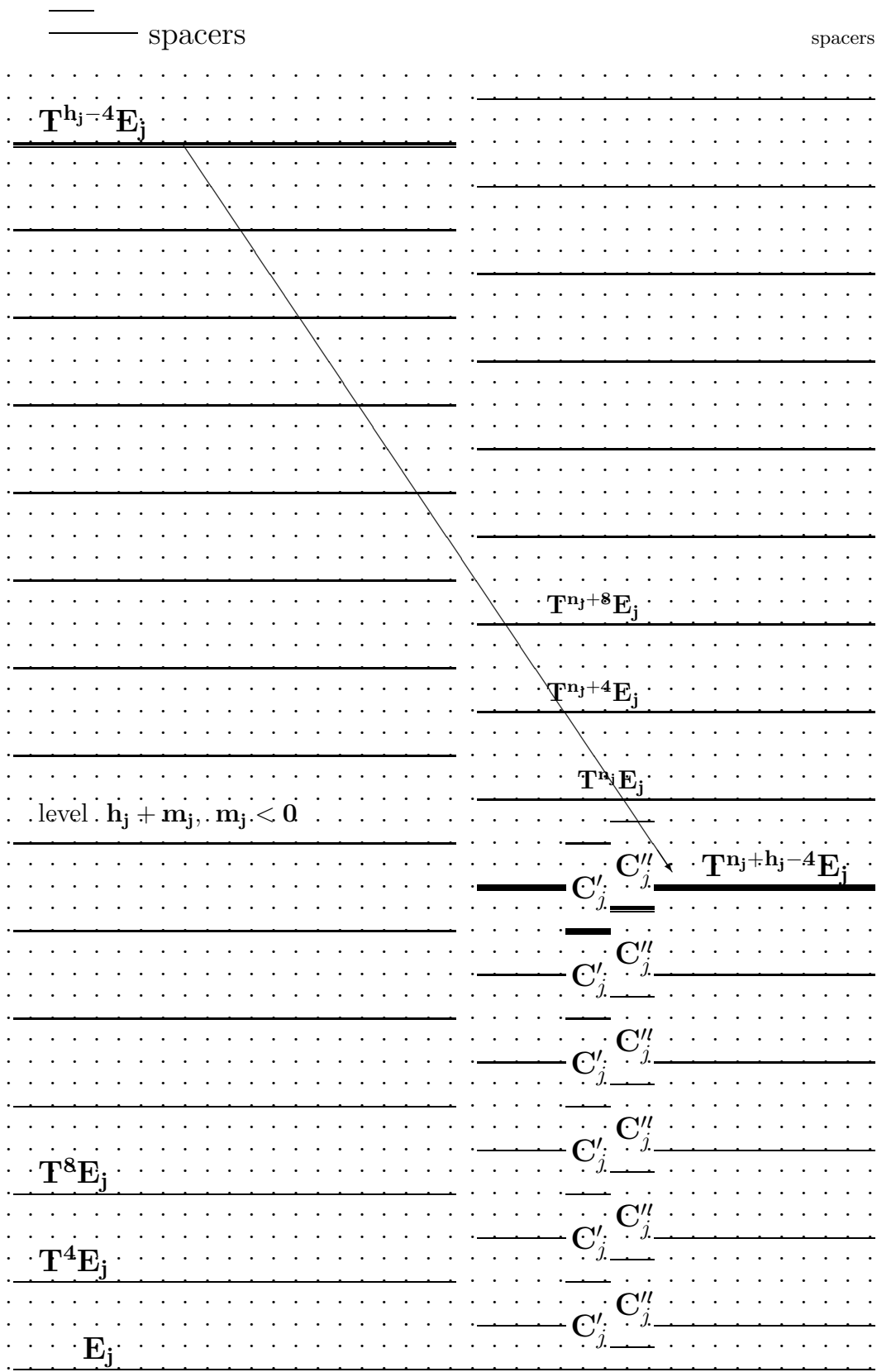
\begin{picture}(0, 630)
\put(0,-20){\it \large  If the  image  at the top of the tower is flat, we are looking for $C_j', C_j''$ at the bottom  }
%\put(22,0){\vector(1,1){170}}
\put(80,560){\vector(2,-3){222}}
\put(0,0){  \line(1, 0){390}}
%\put(-4,600){  \line(1, 0){390}}
\put(16,610){  \line(1, 0){40}}
\put(16,620){  \line(1, 0){20}}
{\Large
\put(65,605){spacers}}

{\small
\put(365,605){spacers}}
%\put(0,0){  \line(0, 1){600}}
%\put(400,0){  \line(0,1){600}}

\multiput(0,0)(10,0){40}{
\multiput(0,0)(0,10){60}%
{${. }$}}

%\multiput(400,0)(0,10){60}%
%{${\bf .}$}
\multiput(0,0)(0,40){15}
{  \line(1,0){200}}

\multiput(210,260)(0,40){9}
{  \line(1,0){180}}

%\multiput(270,590)(0,-40){7}
%{  \line(1,0){20}\put(-15,-25){{\Large \bf C$''_j$}}}
%\multiput(250,560)(0,-40){6}
%{  \line(1,0){20}\put(-15,-25){{\Large \bf C$'_j$}}}

\multiput(210,260)(0,-40){7}
{  \line(1,0){40}}
\multiput(270,250)(0,-40){7}
{  \line(1,0){20}}
\multiput(250,240)(0,-40){7}
{  \line(1,0){20}}
\multiput(250,240)(0,-40){6}
{  \put(6,-25){{\Large \bf C$'_j$}}}
\multiput(250,240)(0,-40){6}
{  \put(26,-15){{\Large \bf C$''_j$}}}
\multiput(290,260)(0,-40){7}
{  \line(1,0){100}}

\put(0,559){  \line(1, 0){200}}
\put(0,558){  \line(1, 0){200}}

\put(210,221){  \line(1, 0){40}}
\put(210,219){  \line(1, 0){40}}
\put(290,221){  \line(1, 0){100}}
\put(290,219){  \line(1, 0){100}}

\put(250,201){  \line(1, 0){20}}
\put(250,199){  \line(1, 0){20}}

\put(270,211){  \line(1, 0){20}}
\put(270,210){  \line(1, 0){20}}
\put(270,209){  \line(1, 0){20}}
%\put(130,100){  \li,ne(1,-2){50}}
{\Large
\put(35,5){${\bf  E_j}$}
		\put(15,45){${\bf  T^4E_j}$}
\put(15,85){${\bf T^8 E_j}$}
}

{\large
\put(255,265){$ \ {\bf T^{n_j} E_j}$}

		\put(245,305){${\bf T^{n_j+4} E_j}$}\put(245,345){${\bf  T^{n_j+8}E_j}$}

}

{\Large
\put(15,565){${\bf  T^{h_j -4}E_j}$}
		\put(315,225){${\bf  T^{n_j+h_j -4}E_j}$}
%\put(150,-12){$r_j=19$}
 }
\put(15,250){\large level $ \ {\bf {h_j+m_j, \ m_j< 0}}$}
\end{picture}
\\ \\ \\ \\ \\
If  the mentioned image of  $E_j$  is close to the top, or to the  bottom,  we use Lemma 1. This finish the proof of Theorem 1.
\newpage

}
\Large
\section{ Bounded constructions and disjointness of powers}
In  \cite{B} J. Bourgain  proved 
 	that bounded constructions  satisfied the Moebius orthogonality property.  This property is a consequence of MSJ  (in fact, of the disjointness  of $T^q$ and $T^p$, $q\neq p$). 

\bf THEOREM 2. \it   If a bounded construction $T$ is weakly mixing and $q\neq p$,  then 
$T^q$ and $T^p$ are disjoint.\rm

Proof.  
 If  roofs over last columns are asymptotically  non-flat, a weakly mixing construction is not rigid, hence,  it has MSJ (Theorem 1).
If   a flatness appears, this means that  our spacer sequence is following: for some sequence of integer intervals  $[\alpha_t, \beta_t]$  $( \beta_t-\alpha_t) \to\infty$

$$ \bar s_j = (s_t, s_t,\dots, s_t, s_t, 0), \ \  \forall  j\in [\alpha_t, \beta_t]      
 $$ 
( "0"-s are necessary in this situation).
From the weak mixing property,  from time to time we see  a break in  the mentioned  behavior (our construction is not an odometer!).
 Thus, infinitely many times  we meet  non-flat roofs over last columns.
Given $\eps>0$, it is not hard  to get for  all $k<\eps^{-1}$  the following weak limits: 
$$T^{kn_j}\to  (1-k\eps) I + k\eps P,$$
where $P=\sum_{i\geq 0}^{\infty} a_i T^i,$   
$\sum_{i\geq 0}^{\infty} a_i =1$, and $a_m> 0$   for some $m>0$.  

If $$T^qJ=JT^p$$ for an operator $J:L_2(\mu)\to L_2(\mu)$, then
$$((1-q\eps) I + q\eps P)J=J((1-p\eps) I + p\eps P).\eqno (\ast)$$

Assume  $J$ to be an indecomposable Markov operator (indecomposable
in the convex set of Markov operators intertwining $T^q$ and $T^p$, or,
in other words,  the corresponding joining is ergodic). All operators $T^nJ,  JT^m$ are indecomposable as well.
Let  $q<p$, from $(\ast)$ it follows $J=JT^{m'} $  for some $m'>0$ (again a  marginal invariance of a joining!).    The ergodicity of 
$T^{m'}$  implies $J=\Theta$.  Thus,  $T^q$ and $T^p$ are disjoint.
 
\large
%%%%%%%%%%%%%%%%%%%%%%%%\put(205,305){$the %\ image \ {\bf T^{n_j} E_j, {\it \ at\ the\ level} \ \  T^{m_j}E_j}$}
\section{Explicit Weak Closure of Actions}
\bf Stochastic Chacon's map. \rm 
Fixing $h_1$,  cut-numbers $r_j\to\infty$,   consider  the set 
of all constructions with spacer sequences
  $s_j(i)\in \{ 0, 1\}$. Let's equip this ensemble
  by Bernoulli measure of type 
$ (a, 1-a)$.
I.e., in a random fashion, with probability $1-a$
 we stack
one spacer  over a column.  Denote   $P=aI+(1-a)T^{-1}$.

\bf THEOREM 3. \it
For  stochastic Chacon's maps $T$  almost surely 
$$ Lim(T)=\{\Theta, P^{m}{P^\ast}^{n}T^k \ : m,n\in\N,\ k\in \Z, m+n\neq 0\}.$$ 

\bf COROLLARY.  \it $ Lim(T)=  \{ \Theta,   T^kP^n : n=0,1,2 ....   n >0, k \in \Z\}$ as $a=0.5$.
\rm

\bf LEMMA. \it 
Let  $T^{n_j} \to Q$,                   $ n_j >0$. 
Represent
$n_j =  q_j h_{p(j)}   +      m_j$ , 
    where    $h_{p(j)}\leq  n_j < h_{p(j)+1}$.    $ 0\leq  m_j \leq  h_{p(j)}$,     $q_j \leq  r_{p(j)}$.
   If $q_j\to\infty $  and  $(r_{p(j)} -q_j)\to  \infty$,  then $Q=\Theta$.
\\
\rm
The proof of this lemma uses the following facts:
For almost all T
for any $q>0$
    $$T^{qh_j} \to  P^q $$ 
and
$$P^q \to \Theta,  \ \ as \ \ q\to \infty,$$ moreover,
  $T^n P^q \to \Theta$  uniformly with respect to $n$.

Assume that $q_j$ (or $r_p(j)- q_j$) to be  bounded.
If   $\eps h_{p(j)}< | m_j |<(1-\eps)h_{p(j)}$, then $Q=\Theta$.  
Let $$\frac {m_j}   {h_p(j)}  \to   0.$$
For all $i$, $0<i < r_{p(j)}-q+1$
$$T^{m_j  -s_j(i)-s_j(i+1)-\dots- s_j(i+q-1)} \approx_w Q,$$
$$T^{n_j}\approx_w \frac 1 {r_{p(j)}} \sum_{i=1}^{r_{p(j)}-q} T^{ -s_j(i)-s_j(i+1)-\dots- s_j(i+q-1)}T^{m_j } \approx_w P^q T^{m_j}.$$

Thus,  $Q=       P^q   \lim_j   T^{m_j}$, then we get 
 $$ Q= P^k {P^*}^n Q^{''}$$
and so on. Note that ${P^*}$ appears in case  of  
$\frac {m_j}   {h_p(j)}  \to   1.$
Iterating, if we cannot stop  by   $ Q^{'''\dots '''}=T^m$, then  $P=\Theta$  ($P^k {P^*}^nT^m \to \Theta$ as $k+n\to \infty$  uniformly with respect to $m$).

\vspace{5mm}
\bf  Bounded staircase flow $T_t$. \rm
Fix $h_1\in R^+$ and cut-sequence $r_j\to\infty$. Define  spacers  by
  $s_j(i)= \frac  {i-1} {r_j}$  ($1\leq i\leq r_j$).

We have 
$$T_{mh_j}\to
 P_m:=\int_{-m}^0T_t dt=T_{-m} P^*_m.$$

\bf THEOREM. \it  For bounded staircase flow $T_t$  
$$Lim(T_t)=
\{\Theta, T_a\prod_{m\in M} P_{m}\ :a\in \R, \ M\subset \N, \ |M|<\infty
 \}.$$ 
\rm

\bf LEMMA. \it 
Let  $T^{t_j} \to Q$,                   $h_{p(j)}\leq  t_j < h_{p(j)+1}$.   
Represent
$t_j =  q_j h_{p(j)}   +      m_j$ , 
    where  $ 0\leq  m_j \leq  h_{p(j)}$,     $q_j \leq  r_{p(j)}$.
   If $q_j\to\infty $  and  $(r_{p(j)} -q_j)\to  \infty$,  then $Q=\Theta$.
\rm

The  proof of the  lemma  is an exercise.

Let  $T_{t_j}\to Q\neq \Theta$,  then Lemma
asserts that there is $q$  such that $t_j =  q h_{p(j)}   +      m_j$.
Again we have 
$$\frac {m_j}   {h_{p(j)}}  \to   0\ \ or \ \ \frac {h_{p(j)}-m_j}   {h_{p(j)}}  \to 1,$$
and
$$ Q= P_q \lim T_{m_j} \ \ \left( \ \ or \ \ Q= {P}_{q+1} \lim T_{m_j-h_{p(j)}}\right).$$
Iterating, taking into account  the fact that 
$T_m\prod_{q\in M} P_{q} \to\Theta $ uniformly with respect to $m\in \R$ as $|M|\to \infty$, we stop as $Q= T_a\prod_{m\in M} P_{m}$.

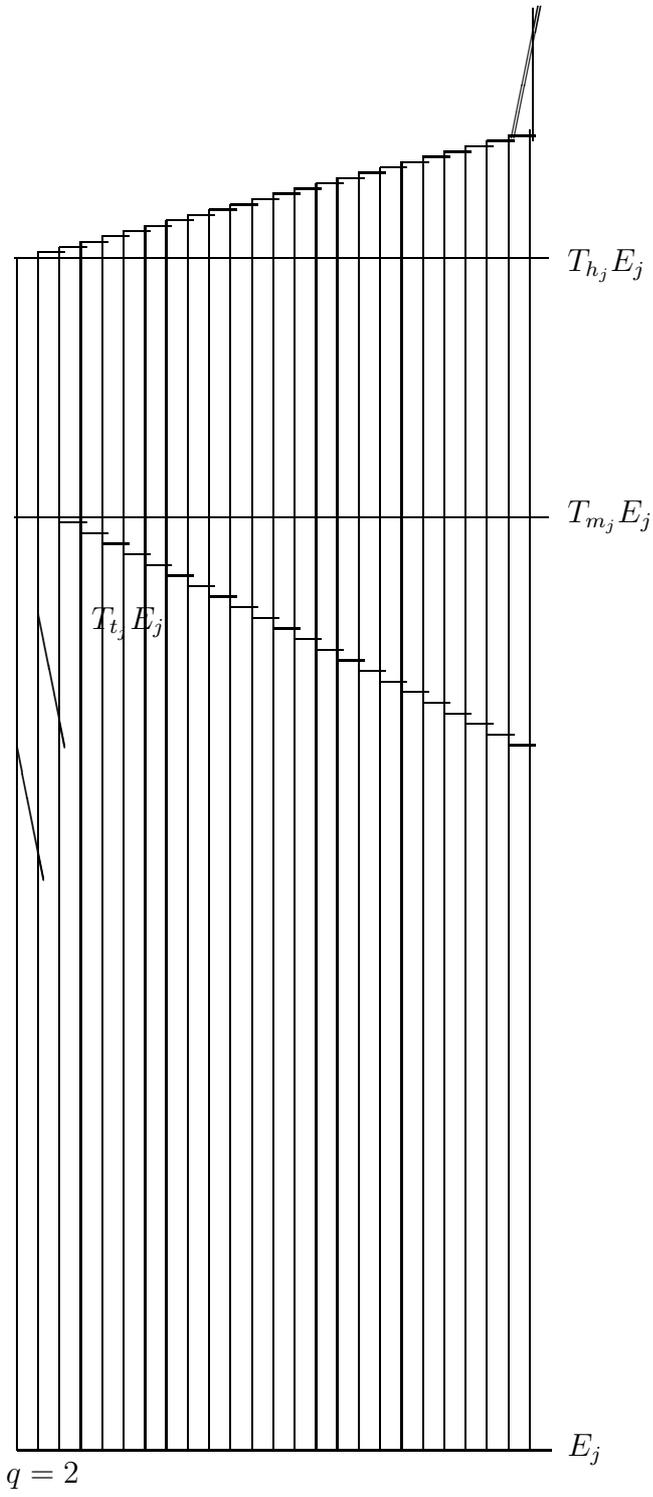
\begin{picture}(650, 490) 
\multiput(0,0)(8,2){25}%
{  \line(0, 1){400}   }
\multiput(0,400)(8,2){24}%%%%%%%%%%%%
{  \line(1, 0){10}   }
\multiput(0,-50)(8,0){25}%
{  \line(0, 1){200}   }
\put(4,-50){\line(1,0){200} }
%\put(2,-49){\line(1,0){198} }
%\put(5,-9){1}
%\put(13,-9){2}\put(21,-9){3}
\multiput(16,300)(8,-4){22}% %%%%%%%%%%%%%
{  \line(1, 0){10}   }
\put(4,215){\line(1,-5){10}}
\put(4,216){\line(1,-5){10}}

\put(12,265){\line(1,-5){10}}
\put(12,266){\line(1,-5){10}}
\put(190,445){\line(1,5){10}}
\put(189,445){\line(1,5){10}}
\put(197,444){\line(0,1){50}}
%\multiput(112,492)(8,-24){10}%
%{  \line(1, 0){10}   }
%\put(38,290){1}
%\put(46,258){2}\put(54,226){3}
%{\large 
%\put(228,140){What is a delay?  If $T^m E^i_j$ is situated in  $T^tE_j$, then }
%\put(228,125){$T^m E^{i+1}_j$ will be  in  $T^{t-d}E_j$. }
%\put(224,110){ This $d$ is an ``upper'' relative  delay}
\put(28,260){  $ T_{t_j}E_j$}
\put(210,300){$T_{m_j}E_j$}
\put(210,395){$T_{h_j}E_j$}
\put(210,-52){$E_j$}
\put(0,-62){$q=2$}
%\put(0,-65){  $E_j^1$, $E_j^2$,$E_j^3$,$\dots$}
%\put(39,15){ delay = 5 }\put(119,145){ delay = 4 }\put(140,445){ delay = 3 }
%\put(250,220){$m =\sum_{i=1}^{d_j} ( h_{j}+i-1) + t_j$
%  } 

\put(3,302){\line(1,0){200} }   \put(3,400){\line(1,0){200} } 
%\put(35,301){\line(1,0){8} }  

\end{picture}  
\Large
\newpage
\section{ Related Problems}

\
1,2.  Thouvenot's questions:   Do    mildly mixing rank-one transformations
possess the properties $PID$? $MSJ(2)$?
\vspace{4mm} \\
3. Is
 $Rank(T^n)$ for  mildly mixing rank-one transformation $T$ equal to $n$ ($n>1$)?   
\vspace{3mm} \\
4. King's  question:  will  any ergodic self-joining of a rank-one transformation be  a limit of off-diagonals measures?
\vspace{3mm} \\
5. Is it true that for any mildly mixing bounded construction $T$ its
symmetric powers $T^{\odot n}$ have simple spectrum? (see \cite{R})
\vspace{4mm} \\
6,7. Let $T$ be  Chacon's map.  Will  the product
$T\otimes T^{2}\otimes T^{3}
\otimes T^{4}\dots $ be of
 simple spectrum  for Chacon's map ? for a mildly mixing  bounded construction? 
\vspace{2mm} \\
8. Consider a rank-one flow called rigid  staircase flow. 
Fix $h_1\in R^+$ and a cut-sequence $r_j\to\infty$. Define  spacers  by
  $$s_j(i)= \frac  {i-1} {jr_j}, \ \   i\leq r_j.$$ Will the corresponding rank-one flow be  simple ?  (see  \cite{JR1} or \cite{T} for definitions  of simplicity. An example of a rigid, simple transformation is given in  \cite{JR}.)  
\vspace{2mm}
\\
9. Let $T$, $T'$  be rank-one constructions possessing  the same  spacer sequences, and $h_1\neq h_1'$ for them.  Are  they disjoint ? 
(It is true for mixing $T$, $T'$, see  	arXiv:1109.0671.)
\vspace{2mm}
\\
10. Are there mildly mixing  (bounded) constructions $T$, $T'$ such that $T'$ is  spectrally isomorphic to $T$  but not isomorphic (as measure-preserving map)  to $T$ and $T^{-1}$?  For mixing rank-one flows  this is
possible:  the flows $T_t$ and $T_{\alpha t }$ ($\alpha >1$)  from  
 arXiv:1002.2808  have the same  spectrum (A. A. Prikhodko), but they are disjoint (see  	arXiv:1109.0671).
\vspace{2mm}
\\
11. For rank-one transformations $MSJ(2)=MSJ$.   Recall that MSJ can be defined as  $MSJ(2)\cap PID$, where the  property $PID$ (pairwise independence
determines  the  global independence)  means  
that any pairwise independent self-joining  have to be trivial (a product measure). This property has been introduced in   \cite{JR1} (see also
\cite{T}).   

  Let $\nu$ be a 3-fold self-joining of weakly mixing rank-one map $T$, and 
for all measurable $A,B\subset X$
$$\nu(A\times B\times X)=\nu(X\times A\times B)=\nu(A\times X\times B)=
\mu(A)\mu(B).$$
Is it true that $\nu= \mu\times \mu\times \mu \ ?$
\newpage
12. Does the condition  $T^{m_j}, T^{n_j},T^{m_j-n_j}\to\Theta$  implies   
$$\mu(A \cap T^{m_j}B \cap  T^{n_j}C)\to \mu(A)
\mu(B)\mu(C)$$
 for a mildly mixing rank-one map $T$?

%\normalsize

 E-mail: vryzh@mail.ru
\end{document}